\documentstyle[12pt]{article}
\newcommand{\const}{\mathop{\rm const}\limits}

\newcommand{\Var}{\mathop{\rm Var}\limits}

\newcommand{\Law}{\mathop{\rm Law}\limits}

\newcommand{\Cov}{\mathop{\rm Cov}\limits}

\begin{document}

\begin{center}

\vspace{3mm}

{\bf Problem of Estimation of Fractional Derivative }\par

\vspace{4mm}

{\bf for a  Spectral Function of Gaussian Stationary Processes} \\

\vspace{4mm}

 $ {\bf E.Ostrovsky^a, \ \ L.Sirota^b } $ \\

\vspace{4mm}

$ ^a $ Corresponding Author. Department of Mathematics and computer science, Bar-Ilan University, 84105, Ramat Gan, Israel.\\

E-mail: eugostrovsky@list.ru\\

\vspace{3mm}

$ ^b $  Department of Mathematics and computer science. Bar-Ilan University,
84105, Ramat Gan, Israel.\\

E-mail: sirota3@bezeqint.net \\

\vspace{4mm}
                    {\sc Abstract.}\\

 \end{center}

 \vspace{3mm}

  We study the problem of nonparametric estimation of the fractional derivative of unknown
 spectral function of Gaussian stationary sequence (time series)
 and show that these problems is well posed with the classical speed of convergence
 when the order of derivative is less than 0.5. \par
  We prove also the asymptotical unbiaseness  and normality of offered estimates with optimal speed of convergence.\par
   For the construction of the  confidence region  in some functional norm we establish the Central Limit Theorem in
 correspondent space of continuous functions  for offered estimates.  \par

 \vspace{4mm}

 \par

\vspace{3mm}


\vspace{4mm}

\section{Notations. Statement of problem.}

\vspace{4mm}

 "Fractional derivatives have been around for centuries  but recently they have
found new applications in physics, hydrology and finance", see  \cite{Meerschaert1}. \par
 Another applications: in the theory of Differential Equations are described in \cite{Miller1};
in statistics see in \cite{Adler1}, \cite{Bapna1}, see also \cite{Golubev1}, \cite{Enikeeva1}, \cite{Ostrovsky107}; in the theory of
integral equations etc. see in  the classical monograph \cite{Samko1}.  \par

\vspace{3mm}

 {\bf  We consider here  the problem of the nonparametric  estimation of the
fractional derivative for a  spectral function of Gaussian stationary sequence. } \par

\vspace{3mm}

 We will prove that if the order of the fractional derivative $  \alpha $ is positive and is strictly less than 1/2,  then this problem
is well posed. In particular, the speed of convergence  of offered (asymptotical) unbiased estimate is $ 1/\sqrt{n}, $ as in the case
of estimation of ordinary spectral function $  F(\lambda); $ they are asymptotical normal still in uniform norm. \par
 Our results improve ones in the books and articles \cite{Anderson1}, chapters 7-9; \cite{Bapna1}, \cite{Borla1}, \cite{Hannan1}, chapter 5;
 \cite{Ibragimov1}, \cite{Ostrovsky1}, chapter 5, section 5.13; \cite{Ostrovsky107} etc. \par

  In particular, it is proved in the preprint \cite{Ostrovsky107}, see also \cite{Bapna1}, that the problem of estimation of fractional
derivative  for function of distribution $  F^{(\alpha)}(x),$  under condition $ \alpha < 1/2 $ is well posed with the optimal rate of convergence
$  1/\sqrt{n}, $ when $  n  $ is the volume of the sample, and is announced an analogues fact about the fractional derivative of spectral function. \par

\vspace{3mm}

  Let $ \alpha = \const \in (0,1); $ and let also $ g = g(x), \  x \in R $  be certain measurable numerical function. We recall for reader
convenience  that the fractional derivative of a Rieman-Liouville type of order $  \alpha, \ \alpha \in (0,1): \  D^{\alpha}[g](x)  = g^{(\alpha)}(x) $
is defined as follows: $ \Gamma(1-\alpha)  g^{(\alpha)}(x) =  $

$$
\Gamma(1-\alpha) \ D^{\alpha}[g](x)  =
\Gamma(1-\alpha) \ D^{\alpha}_x[g](x) \stackrel{def}{= } \frac{d}{dx} \int_0^x \frac{g(t) \ dt}{(x-t)^{\alpha}}, \eqno(1.1)
$$
see, e.g. the classical monograph of S.G.Samko, A.A.Kilbas and O.I.Marichev \cite{Samko1}, pp. 33-38; see also \cite{Miller1}.\par

 Hereafter $ \Gamma(\cdot) $ denotes the ordinary $ \Gamma \ $ function. \par

 We agree to take $ D^{\alpha}[g](x_0) = 0,  $ if at the point $  x_0 $ the expression $  D^{\alpha}[g](x_0)  $ does not exists. \par

Notice that the operator of the fractional derivative is non-local, if $ \alpha $ is not integer non-negative number.\par

  Note also that for the considered further functions this fractional derivative there exists almost everywhere. \par

 Recall also that the fractional  integral $  I^{(\alpha)}[\phi](x) = I^{\alpha}[\phi](x)  $ of a Rieman-Liouville type of
an order $ \alpha, 0 < \alpha < 1 $ is defined as follows:

$$
I^{(\alpha)}[\phi](x) \stackrel{def}{=} \frac{1}{\Gamma(\alpha)} \cdot \int_0^x \frac{\phi(t) \ dt}{(x-t)^{1 - \alpha}}, \ x,t > 0. \eqno(1.2)
$$
 It is known  (theorem of Abel, see \cite{Samko1}, chapter 2, section 2.1)
 that the operator $ I^{(\alpha)}[\cdot]  $ is inverse to the fractional derivative operator $ D^{(\alpha)}[\cdot],  $
at least  in the class of absolutely continuous functions. \par

\vspace{3mm}

  Evidently, (Kolmogorov's theorem), the problem of distribution function estimation $ (\alpha = 0)  $ is well posed in the uniform norm.
V.D.Konakov in \cite{Konakov1}  proved in contradiction that  the problem of the spectral density estimation, i.e. when $ \alpha = 1, $
is ill posed. \par
  Roughly speaking, the result of V.D.Konakov  may be reformulated as follows. Certain problem of statistical estimation is well posed
iff there exists an {\it  continuous } in natural distance estimate (more exactly, a sequence  of such continuous estimates) such that
the speed of convergence is equal (or less than) $  1/\sqrt{n}. $ As a rule these estimations are asymptotically normal.\par

\vspace{4mm}

\section{ Main result: estimation of fractional derivatives of spectral function.}

 \vspace{4mm}

 Let us consider in this section the classical problem of estimation of  fractional derivatives for spectral function,
or equally the problem of fractional integral for spectral density estimation. \par

  Let $ \eta_k, \ k = 1,2,\ldots,n $ be real valued, centered: $  {\bf E} \eta(k) = 0 $ Gaussian distributed stationary random
 sequence (time series) with (unknown) even covariation function $  r = r(m), $  spectral function
 $ F(\lambda), \ \lambda \in [0, 2 \pi],  \ F(0+) = F(0) = 0, $ and with  spectral density $  f(\lambda), $ (if there exists):

$$
r(m) = \Cov(\eta(j+m), \eta(j)) =  {\bf E} \eta(j+m)\cdot \eta(j) =
$$

$$
\int_{[-\pi,  \pi]} \cos(\lambda m) dF(\lambda) =
\int_{[-\pi,  \pi]} \cos(\lambda m) f(\lambda) \ d \lambda, \eqno(2.1)
$$
so that

$$
F(\lambda) = \int_0^{\lambda} f(t) dt = I^{(1)}[f](\lambda); \ F(0) = F(0+) = 0.
$$

 We can and will suppose without loss of generality that the spectral density,  i.e. the
function $  f = f(\lambda) $ to be continued on the whole axis $  R^1 $
as a $  2 \pi $ periodical (continuous) function  $  f(\lambda) = f(\lambda + 2 \pi). $ \par

\vspace{3mm}

 The periodogram of this sequence will be denoted by $ J_n(\lambda), \ 0 \le \lambda \le 2 \pi: $

$$
J_n(\lambda) := (2 \pi n)^{-1} \left| \sum_{k=1}^n e^{i k \lambda} \eta(k)  \right|^2; \ i^2 = - 1. \eqno(2.2)
$$

 We intend here to estimate the fractional derivative  $  F^{(\alpha)}(\lambda)  $ of the spectral function  $  F(\lambda). $ \par
 Recall that the problem of $  F(\cdot) $ estimation is well posed, theorem of I.A.Ibragimov, \cite{Ibragimov1},
 in contradiction to the problem of spectral density $ f(\cdot) = F^{(1)}(\cdot) $  estimation.\par
 We assume as before $ 0 < \alpha < 1/2,  $ and denote $  \beta = 1 - \alpha; \ \beta \in (1/2, 1). $  \par

Heuristic arguments. We have using the group properties of the fractional derivative-integral operators

$$
F^{(\alpha)} = D^{\alpha}[F] = D^{\alpha} I^1[f] =D^{\alpha} D^{-1}[f] = D^{\alpha - 1}[f] =
$$

$$
I^{1 - \alpha}[f] = I^{\beta}[f] \approx  I^{\beta} [J_n].\eqno(2.3)
$$

 Thus, we can offer as an estimation of $ F^{(\alpha)} $ the following statistics

$$
F_{\alpha,n}(\lambda) :=  I^{\beta} [J_n](\lambda) = \frac{1}{\Gamma(\beta)}\int_0^{\lambda} \frac{J_n(t) \ dt}{(t - \lambda)^{1-\beta }} =
\frac{1}{\Gamma(1-\alpha)}\int_0^{\lambda} \frac{J_n(t) \ dt}{(t - \lambda)^{\alpha }}. \eqno(2.4)
$$

 Let us introduce a sequence  as $  n = 1,2,\ldots $  of the normed and mean zero random processes

$$
\zeta_n(\lambda) := \sqrt{n} \left(  F_{\alpha,n}(\lambda) - {\bf E} F_{\alpha,n}(\lambda)  \right) \eqno(2.5)
$$
and also a Gaussian centered separable random process $ \zeta(\lambda) = \zeta_{\infty}(\lambda)  $ with covariation function

$$
\Theta^{(f)}_{\alpha}(\lambda, \mu)  =
\Theta_{\alpha}(\lambda, \mu) \stackrel{def}{=} \Cov(\zeta_{\infty}(\lambda), \zeta_{\infty}(\mu)) = {\bf E}\zeta(\lambda) \cdot \zeta(\mu) :=
$$

$$
\frac{4 \pi}{\Gamma^2(1 - \alpha)} \cdot
\int_0^{ \lambda \cap \mu } \frac{f^2(\nu) \ d \nu}{(\lambda - \nu)^{\alpha} \ (\mu - \nu)^{\alpha}}, \eqno(2.6)
$$
so that

$$
\sigma^2_{\alpha}(\lambda) \stackrel{def}{=}
\Var \zeta_{\infty}(\lambda) =  \Theta_{\alpha}(\lambda, \lambda) =
$$

$$
\frac{4 \pi}{\Gamma^2(1 - \alpha)} \cdot \int_0^{\lambda}
\frac{f^2(\nu) d \nu}{ (\lambda - \nu)^{2 \alpha} } =
\frac{4 \pi \Gamma(1 - 2 \alpha)}{\Gamma^2(1 - \alpha)} \cdot I^{2 \alpha}[f^2](\lambda) < \infty, \eqno(2.7)
$$
as long as the function $  f = f(\lambda) $ is presumed to be continuous and $  \alpha < 1/2. $\par

\vspace{4mm}

{\bf Theorem 2.1.} {\it Suppose as before $ 0 < \alpha < 1/2 $ and
that the spectral density  $  f(\lambda) $ there exists and is continuous and strictly positive on the (closed) circle $ [0, 2 \pi]: $

$$
\exists C_1, C_2, \ 0 < C_1 < C_2 < \infty \ \Rightarrow C_1 \le f(\lambda) \le C_2, \eqno(2.8)
$$
in particular

$$
f(0) = f(0+) = f(2 \pi -0)  = f(2 \pi).
$$

 Our statement: the sequence of the distributions of continuous random processes  $  \zeta_n(\cdot) $ converges weakly as
 $ n \to \infty $ in the space of continuous periodical functions $  C^*(0, 2 \pi), $  i.e. in the Prokhorov-Skorokhod sense,
to the distribution of the continuous random process  $  \zeta_{\infty}(\cdot): $  for arbitrary  continuous bounded functional}
 $  G:  C^*(0, 2 \pi)  \to R $

$$
\lim_{n \to \infty} {\bf E} G(\zeta_n(\cdot)) =  {\bf E} G(\zeta_{\infty}(\cdot)). \eqno(2.9)
$$

\vspace{3mm}

{\bf Proof.}\\

\vspace{3mm}

{\bf 1.} Note first of all that all the  considered here r.p. $  \zeta_n(\lambda), \ \zeta_{\infty}(\lambda)  $ are continuous
with probability one. Indeed, the covariation function $ \Theta_{\alpha}(\lambda, \mu)  $ satisfies the H\"older's condition with
exponent $ \alpha - \epsilon  $ for arbitrary value $  \epsilon \in (0, \alpha/2). $ For instance, the function
$ \lambda \to \sigma^2(\lambda) =\sigma^2_{\alpha}(\lambda)  $  in (2.7)  being represented as a fractional integral from the
continuous function $  f^2(\cdot), $ obeys the module continuity of a form

$$
\omega(\sigma^2, h) \le C(\alpha) \ h^{\alpha - 1/p} \ |f^2|_p, \ p > 2/\alpha, \eqno(2.10)
$$

$$
|g(\cdot)|_p := \left[ \int_0^{2 \pi} |g(\lambda)|^p \ d \lambda    \right]^{1/p},
$$
see \cite{Samko1}, chapter 1, section 3.3., pp. 66-71.\par

 The general case $ \lambda \ne \mu $ may be establish analogously using  at the same  statement in the book \cite{Samko1}. \par

 Hereafter as ordinary $ \omega(g(\cdot),h), \ h \in [0, 2 \pi] $ denotes the module of continuity of (possible continuous and periodical)
function $  g = g(\lambda): $

$$
\omega(g(\cdot),h) \stackrel{def}{=} \sup_{ |\lambda - \mu| \le h } |g(\lambda) - g(\mu)|. \eqno(2.11)
$$

\vspace{3mm}

{\bf 2.} I.A.Ibragimov in  proved in  \cite{Ibragimov1} in particular the Central Limit Theorem for the sequence of r.p.

$$
\tau_n(\lambda) := \sqrt{n} \left[ F_n(\lambda) - {\bf E} F_n(\lambda) \right]
$$
in the space $  C^*(0, 2 \pi), $  where

$$
F_n(\lambda) = \int_0^{\lambda} J_n(t) \ dt \eqno(2.12)
$$
is  ordinary empirical spectral function. \par

 To be more precise, we introduce following I.A.Ibragimov  also a centered separable Gaussian random process
$  \tau(\lambda), \ \lambda \in [0, 2 \pi] $ with  covariation function

$$
\Cov(\tau(\lambda), \tau(\mu)) = 4 \pi \ \int_0^{\min(\lambda, \mu)} f^2(x) \ dx. \eqno(2.13)
$$

 Define a non - negative function $ \beta = \beta(\lambda), \ \lambda \in [0, 2 \pi] $

$$
\beta^2(\lambda) := 4 \pi \ \int_0^{\lambda} f^2(x) \ dx, \eqno(2.14)
$$
and  let us introduce also the following continuous distance function

$$
d_{\beta}(\lambda, \mu) := |\beta(\lambda) - \beta(\mu)|. \eqno(2.15)
$$

 Then the r.p. $ \tau(\lambda) $ may be represented in the sense of distributional coincidence  as follows:

$$
\tau(\lambda) \stackrel{d}{=} B(\beta(\lambda)),
$$
where $  B(t) $ is ordinary Brownian motion (Wiener's process).\par

\vspace{3mm}

 Obviously, the r.p. $ \tau(\cdot) $ is continuous a.e. \par

 \vspace{3mm}

 I.A.Ibragimov proved also that as $  n \to \infty  $  in the space $ C^*[0, 2 \pi]  $

$$
 \Law(\tau_n(\cdot)) = \Law( \sqrt{n} \left[ F_n(\cdot) - {\bf E}F_n(\cdot) \right] )  \to \Law (\tau(\cdot)). \eqno(2.16)
$$
See also \cite{Buldygin2}, \cite{Dahlhaus1}, \cite{Levit1}.

 Therefore, the finite-dimensional distributions of r.p. $ \zeta_n(\cdot) $ converges to ones for the r.p.  $ \zeta(\cdot), $
which are of course also Gaussian. \par

Moreover,

$$
\sup_n \sup_{\lambda \in [0, 2 \pi]} {\bf E} \tau_n^2(\lambda) \le C_1 < \infty, \eqno(2.17a)
$$

$$
\sup_n {\bf E} \left[  \tau_n(\lambda) - \tau_n(\mu)  \right]^2  \le C_2 \  d_{\beta}(\lambda, \mu), \eqno(2.17b)
$$
see  \cite{Ibragimov1}. \par

\vspace{4mm}

{\bf 3.} Let us calculate the covariation function of the r.p. $  \tau = \tau(\lambda). $ \par

 More detail, we propose

$$
\lim_{n \to \infty} n \cdot \Cov \left\{ F_{\alpha,n}(\lambda), F_{\alpha,n}(\mu)  \right\} =
\frac{4 \pi \Gamma(1 - 2 \alpha)}{\Gamma^2(1 - \alpha)} \cdot D^{-\alpha}_{\lambda} \left[ D^{-\alpha}_{\mu}[f^2]  \right]=
$$

$$
\frac{4 \pi}{\Gamma^2(1 - \alpha)} \cdot
\int_0^{ \lambda \cap \mu } \frac{f^2(\nu) \ d \nu}{(\lambda - \nu)^{\alpha} \ (\mu - \nu)^{\alpha}}
=: \Theta_{\alpha}(\lambda, \mu). \eqno(2.18)
$$

Note that the last integral is finite since the function $ f $ is bounded and  $  \alpha < 1/2. $ \par

\vspace{3mm}

 This  assertion (2.18) follows immediately from
 the following proposition, see  the  fundamental monograph of T.W.Anderson
\cite{Anderson1}, chapter 5, page 564-572, theorem 9.3.1:
if $ w(\lambda, \nu) $ is non-negative integrable relative the second variable function, then

$$
\lim_{n \to \infty} {\bf E} \int_0^{2 \pi} w(\lambda, \nu) J_n(\nu) d \nu =
\lim_{n \to \infty}  \int_0^{2 \pi} w(\lambda, \nu) {\bf E} J_n(\nu) d \nu =
$$

$$
\int_0^{2 \pi} w(\lambda, \nu) f(\nu) d \nu,
$$

$$
\lim_{n \to \infty} n \cdot \Cov \left( \int_0^{2 \pi} w(\lambda_1, \nu) J_n(\nu) d \nu,
\int_0^{2 \pi} w(\lambda_2, \nu) J_n(\nu) d \nu  \right) =
$$

$$
4 \pi \int_0^{2 \pi} w(\lambda_1, \nu) \ w(\lambda_2, \nu) \ f^2(\lambda) \ d \lambda,
$$
with remainder terms. We choose $ w(\lambda, \nu) = |\lambda - \nu|^{-\alpha}; $ it is easy to verify that
all the conditions of the mentioned result are satisfied. \par
 Recall also that the considered stationary sequence $  \{ \eta(k) \}  $ is Gaussian, i.e. without cumulant function. \par

\vspace{3mm}

{\bf 4.} It remains only to establish the weak compactness of the distributions $ \zeta_n(\cdot)  $  in the space of continuous
functions $  C^*[0, 2 \pi]. $  Note first of all that

$$
C_3 |\lambda - \mu| \le d_{\beta}(\lambda, \mu) \le C_4 |\lambda - \mu|. \eqno(2.19)
$$
 Further,  the r.p. $ \tau_n(\lambda) $ can be represented as a two-dimensional mean zero stochastic integral over a Gaussian
stochastic measure,  or equally  in the terminology of the article  \cite{Kozachenko2} "square Gaussian random  vectors (variables)" (SGV)
 or more generally "square Gaussian random  process" (SGP). \par

{\bf 5.} We will use some facts about these random  vectors and processes, see \cite{Buldygin1}, \cite{Kozachenko2}.  Introduce
following V.V.Buldygin and Yu.V.Kozachenko the next function

$$
\phi(\lambda) := \ln \left[ (1 - |\lambda|)^{-1/2} \ e^{-|\lambda|/2}  \right], \ |\lambda| < 1, \eqno(2.20)
$$
 and $  \phi(\lambda)  = + \infty  $ otherwise. V.V.Buldygin and Yu.V.Kozachenko in \cite{Buldygin1} have proved that

$$
\sup_{\eta \in SGV }  {\bf E} \exp \left( \lambda \eta/\sqrt{2 \Var \eta} \right) \le \exp{\phi(\lambda)}, \eqno(2.21)
$$
where in (2.21) the supremum is calculated over all the non-trivial random variable $ \eta $ from   the set SGV,
which may be defined  on arbitrary probability space. \par

 Since for the r.v. of the form $  \eta = \pm (\xi^2 - 1),  $ where the r.v. $  \xi   $ has a standard normal distribution, in the
inequality (2.21) take place the equality, we conclude

$$
\sup_{\eta \in SGV }  {\bf E} \exp \left( \lambda \eta/\sqrt{2 \Var \eta} \right) = \exp{\phi(\lambda)}. \eqno(2.21a)
$$

{\bf 6.} The relations (2.21) (and (2.21a)) may be transformed.  The function $  \phi = \phi(\lambda) $ generated so - called
Banach space $  B(\phi)  $ as follows.\par

   We say by definition that the centered random variable (r.v) $  \eta $ defined on some sufficiently rich probability space
 belongs to the space $  B(\phi), $ if there exists some non-negative constant $ \gamma \ge 0 $ such that

$$
\forall \lambda \in R \ \Rightarrow {\bf E} \exp(\lambda \ \eta) \le \exp(\phi(\gamma \ \lambda)). \eqno(2.22)
$$

 The minimal value $  \gamma $ satisfying (2.22) is called a $   B(\phi) $ norm of the variable $ \eta, $  write

 $$
 ||\xi||B(\phi) = \inf \{ \gamma, \ \gamma > 0: \ \forall \lambda \ \Rightarrow
 {\bf E}\exp(\lambda \xi) \le \exp(\phi(\lambda \ \gamma)) \}. \eqno(2.23)
 $$

  The space $ B(\phi) $ with respect to the norm $ || \cdot ||B(\phi) $ and
ordinary algebraic operations is a Banach space which is isomorphic to the subspace
consisted on all the centered variables from the so-called {\it  exponential} Orlicz’s space
$ (\Omega,F,{\bf P}), N(\cdot) $ with $ N \ - $ function

$$
N(u) = \exp(\phi^*(u)) - 1, \ \phi^*(u) := \sup_{\lambda} (\lambda u - \phi(\lambda)).
$$
The detail investigation of alike spaces see in  \cite{Kozachenko1}, \cite{Ostrovsky1}, \cite{Ostrovsky8}. \par

\vspace{3mm}

{\bf 7.} The Ibragimov's results (2.17a)  and   (2.17b) may be rewritten as follows

$$
\sup_n \sup_{\lambda \in [0, 2 \pi]} || \tau_n(\lambda)|| B(\phi) \le C_5 < \infty, \eqno(2.24a)
$$

$$
\sup_n  ||  \tau_n(\lambda) - \tau_n(\mu) ||B(\phi)  \le C_6 \  d^{1/2}_{\beta}(\lambda, \mu), \eqno(2.24b)
$$
as long as $  \tau_n(\lambda) $ is a two - dimensional stochastic integral. \par

 Recall, see \cite{Kozachenko1}, \cite{Ostrovsky1} that the norm $  || \cdot||B\phi  $ is equivalent on the (closed) subspace
of mean  zero random variables with the following norm

$$
|||\eta|||G\psi \stackrel{def}{=} \sup_{p \ge 2} \left[ \frac{|\eta|_p}{\psi(p)} \right], \eqno(2.25)
$$
where $ \psi(p) = p.  $ This means that on the set of centered variables from the space  $ B(\phi) \  S^o:=
\{ \eta: {\bf E} \eta = 0  \} \cap B(\phi) $

$$
 0 < \inf_{\eta \in S^o} \left[ \frac{||\eta||B\phi}{|||\eta|||G\psi} \right] \le
  \sup_{\eta \in S^o} \left[ \frac{||\eta||B\phi}{|||\eta|||G\psi} \right]   < \infty. \eqno(2.26)
$$

  Therefore

$$
\sup_n \sup_{\lambda \in [0, 2 \pi]} || \tau_n(\lambda)|| G\psi \le C_7 < \infty, \eqno(2.27a)
$$

$$
\sup_n  ||  \tau_n(\lambda) - \tau_n(\mu) ||G\psi  \le C_8 \  d^{1/2}_{\beta}(\lambda, \mu) \eqno(2.27b)
$$
 for some  positive finite constants $  C_7, \ C_8.  $ \par

 \vspace{3mm}

{\bf 8.}  Let now $  p  $  be fixed number  greatest than $  1/\alpha; $ for example

$$
p := p_0 \stackrel{def}{=} \frac{8}{1 - 2 \alpha}. \eqno(2.28)
$$
  We deduce on the basis of inequalities (2.27a) and (2.27b)

$$
\sup_n \sup_{\lambda \in [0, 2 \pi]} | \tau_n(\lambda)|_p \le C_9 p < \infty, \eqno(2.29a)
$$

$$
\sup_n  |  \tau_n(\lambda) - \tau_n(\mu) |_p   \le C_{10} \  p \ |\lambda -  \mu|^{1/2}. \eqno(2.29b)
$$

\vspace{3mm}

 We intend to apply the so-called {\it majorizing measures method,} see e.g.  \cite{Ostrovsky109}, choosing
as a capacity of the  majorizing measure the ordinary Lebesgue measure. The direct application of the proposition 2.1
from \cite{Ostrovsky109} gives us in the considered case the estimation

$$
|\tau_n(\lambda) - \tau_n(\mu) | \le C(p) \ X(n,p) \ |\lambda - \mu|^{1/2 - 2/p}, \eqno(2.30)
$$
or equally

$$
\omega(\tau_n(\cdot), h)  \le C(p) \ X(n,p) \ h^{1/2 - 2/p}, \ 0 \le h \le 2 \pi, \eqno(2.30a)
$$
where  the sequence as $  n = 1,2,\ldots  $ of non - negative r.v. $ X_{n,p}  $ is such that

$$
\sup_n {\bf E} X^p_{n,p}  = 1. \eqno(2.31)
$$

 Note that $  p_0 = 8/(1 - 2 \alpha)  $  and a fortiori $ 1/2 - 2/p_0 > 0. $\par

 \vspace{3mm}

{\bf 9.} We know that the r.p. $ \zeta_n(\cdot)  $ is the fractional derivative  from $ \tau_n(\cdot):$

$$
\zeta_n(\lambda) = D^{\alpha} \tau_n(\lambda)
$$
and in addition $  \tau_n(0) = 0.  $ We apply the inequality for such a functions

$$
\omega(D^{\alpha}f, h) \le C(\alpha) \cdot \int_0^h \frac{\omega(f,t) \ dt}{t^{1 + \alpha} }, \eqno(2.32)
$$
see \cite{Samko1}, p. 250-253, theorem 3.16;  and get

$$
\omega(\zeta_n, h) \le C \ X_{n,p} \ \int_0^h \frac{t^{1/2 - 2/p} \ dt  }{t^{ 1 + \alpha  }} \le
C(\alpha) \ X_{n,p} \ h^{(1-2\alpha)/8}. \eqno(2.33)
$$
 Since $ \alpha \in (0,1/2), $ we conclude taking into account (2.31) that the sequence of r.p. $ \zeta_n(\cdot)  $ satisfies the famous
Prokhorov's criterion  \cite{Prokhorov1} for weak compactness of the (Borelian) probability measures in the space of
continuous functions.\par

 This completes the proof of theorem 2.1.\par

 \vspace{4mm}

 {\bf Theorem  2.2.} {\it Let all the conditions of theorem 2.1. be satisfied.  Suppose in addition

 $$
 \lim_{n \to \infty} [ \sqrt{n} \ \omega(f, 1/n) \ |\ln \omega(f, 1/n)| ] = 0. \eqno(2.34)
 $$

 We propose that the sequence of the distributions of continuous random processes

 $$
  \theta_n(\lambda) \stackrel{def}{=} \sqrt{n} \ (F_{n,\alpha}(\lambda) - F^{(\alpha)}(\lambda))
  $$
 converges weakly as $ n \to \infty $ in the space of continuous periodical functions $  C^*(0, 2 \pi), $
i.e. in the Prokhorov-Skorokhod sense, to the distribution of at the same continuous random process  }
$  \zeta_{\infty}(\cdot). $   \par

\vspace{3mm}

{\bf Proof.} It is sufficient to justify that

$$
\lim_{n \to \infty} \sup_{\lambda \in [0, 2 \pi]} | {\bf E} J_n(\lambda)  -  F(\lambda)| = 0. \eqno(2.35)
$$

The expression for $ {\bf E} J_n(\lambda) $ is given and investigated, e.g.,  in \cite{Anderson1}, chapter 8, sections 8.2 - 8.3:

$$
{\bf E} J_n(\lambda) = \int_{-\pi}^{\pi} \Phi_n(\lambda - \nu) \ f(\nu) \ d \nu = [\Phi_n * f](\lambda), \eqno(2.36)
$$
where $ \Phi_n(\cdot)  $ is the well - known Fejer's kernel

$$
\Phi_n(\lambda) = \frac{\sin^2(n\lambda/2) }{2 \pi n \ \sin^2(\lambda/2)}.
$$
 The error of the Fejer's approximation $ \Phi_n * f  - f $ in the uniform norm is investigated in many works: \cite{Burkill1},
\cite{Nishishiraho1},  \cite{Vore1}, pp. 339 - 341,  \cite{Vore2} etc.  For instance,

$$
\sup_{\lambda \in [0, 2 \pi]} |[\Phi_n * f](\lambda)  - f(\lambda)| \le C \ \omega(f, 1/n) \ |\ln \omega(f, 1/n)|,
$$
see \cite{Burkill1}, pp. 33 - 40. Therefore, the equality (2.35) there holds by virtue of condition (2.34). \par

 \vspace{3mm}

Note that the condition (2.34) is satisfied if for example

$$
\omega(f, h) \le C \ h^{\beta}, \ 0 \le h \le 2 \pi, \ \beta = \const > 1/2, \eqno(2.37).
$$
see \cite{Burkill1}, pp. 33 - 40. Note only that the condition (2.37)  is very weak.\par
 T.W.Anderson in \cite{Anderson1}, chapter 8, section 8.3 imposed on the  spectral density  $  f(\lambda) $ a more
strong restriction

$$
\sum_{m=0}^{\infty}  \sigma(m) < \infty,
$$
 where

$$
f(\lambda) = \sum_{m=0}^{\infty}  \sigma(m) \ \cos(m \lambda).
$$

\vspace{4mm}

{\bf Remark 2.1.} Emerging in the equality (2.7) the variable

$$
 I^{2 \alpha}[f^2](\lambda) = \frac{1}{\Gamma(2 \alpha)} \int_0^{\lambda} \frac{f^2(\nu) d \nu}{(\lambda - \nu)^{1 - 2 \alpha}}
$$
may be $ n^{-1/2} \ - $ consistent estimated as follows:

$$
 I^{2 \alpha}[f^2](\lambda) \approx \frac{1}{\Gamma(2 \alpha)} \int_0^{\lambda} \frac{J_n^2(\nu) d \nu}{(\lambda - \nu)^{1 - 2 \alpha}}.
$$

\vspace{3mm}

{\bf Remark 2.2.} We have proved the asymptotical normality under certain condditions of the sequence of random processes

$$
\theta_n(\lambda) = \sqrt{n} \left\{ F_{n,\alpha}(\lambda) - F^{(\alpha)}(\lambda) \right\}
$$
 as $  n \to \infty $ in the space $ C^*(0, 2 \pi) $ of continuous
functions.  Therefore if the value $ n $ is "sufficiently great"

$$
{\bf P} \left(  \sqrt{n} \cdot \max_{\lambda} \left| \left\{ F_{\alpha,n}(\lambda) - F^{(\alpha)}(\lambda) \right\} \right|  > u  \right)
\approx {\bf P} (\max_{\lambda} | \zeta_{\infty}(\lambda) | > u), \ u = \const > 0.
$$

 The asymptotical as $  u \to \infty  $ behavior of the last probability is fundamental investigated in the monograph
\cite{Piterbarg2}, see also  \cite{Piterbarg1}:

$$
 {\bf P} (\max_{\lambda} | \zeta_{\infty}(\lambda) | > u) \sim H(\alpha) \ u^{\kappa - 1} \ \exp \left(-u^2/\sigma^2 \right),
$$

$$
H(\alpha), \ \kappa = \const, \ \sigma^2 = \sigma^2(\alpha) = \max_{\lambda \in (0, 2 \pi)} \Theta_{\alpha}(\lambda,\lambda).
$$

 The last equalities may be used by construction of confidence region  for $ F^{(\alpha)}(\cdot) $ in the uniform norm.
Indeed, let $ 1 - \delta $  be the reliability of confidence region, for example, $ 0.95  $ or $ 0.99 $ etc.
Let $ u_0 = u_0(\delta) $ be a maximal root of the equation

$$
H(\alpha) \ u_0^{\kappa - 1} \ \exp \left(-u_0^2/\sigma^2 \right) = \delta,
$$
then with probability $ \approx 1 - \delta $

$$
 \sup_{\lambda \in (0, 2 \pi)} \left| F_{\alpha,n}(\lambda) - F^{(\alpha)}(\lambda)  \right| \le \frac{u_0(\delta)}{\sqrt{n}}.
$$

\vspace{4mm}

 \section{CLT in H\"older spaces for spectral functions estimation.}

\vspace{3mm}

 Let $ (X = \{x \},d) $ be compact metric space relative some distance (or semi-distance) $  d = d(x_1, x_2). $
 The {\it modified}  H\"older (Lipshitz) space  $ H^o(d)  $  consists by definition on all the numerical (real or complex)
 continuous relative the distance $ d = d(x_1,x_2) $ functions $ f: X \to R $  satisfying the addition condition

$$
\lim_{\delta \to 0+} \frac{\omega(f,d, \delta)}{\delta} = 0.  \eqno(3.1)
$$
 Here $ \omega(f, \delta) = \omega(f,d, \delta)  $  is  as before uniform module of continuity of the (continuous) function $  f $
relative the distance (metric) $ d(\cdot, \cdot): $

$$
\omega(f,d, \delta) = \omega(f,\delta) = \sup_{x_1,x_2: d(x_1,x_2) \le \delta} |f(x_1) - f(x_2)|. \eqno(3.2)
$$

The norm of the space $ H^o(d) $  is defined as follows:

$$
||f||H^o(d) = \sup_{x \in X} |f(x)| + \sup_{d(x_1, x_2) > 0} \left\{ \frac{|f(x_1) - f(x_2)|}{d(x_1, x_2)}  \right\}. \eqno(3.3)
$$

   The detail investigation of these spaces with applications in the theory of non-linear singular integral equations
is undergoing in the first chapter of a monograph  of  Gusejnov A.I., Muchtarov Ch.Sh.
\cite{Gusejnov1}. We itemize some used facts about these spaces.\par

 This modification of the classical H\"older (Lipshitz) space $ H(d) $ in which the condition (3.1) do not be  presumed
and hence  is not separable, is really {\it separable} Banach space, in particular is  linear, normed and complete and in turn
is a closed subspace of $  H(d). $ \par

 Note but the space $ H^o(d) $  may be trivial, i.e. may consists only on constant functions. Let for instance, $  X  $ be convex
connected closed bounded domain in the space $  R^m, \ m = 1,2, \ldots $ and let $  d(x_1, x_2) = |x_1 - x_2| $ be usual Euclidean distance.
Then the space $ H^o(d)  $ is trivial: $  \dim H^o(d) = 1. $\par
  The space $  H^o(d^{\beta}), \ \beta = \const \in (0,1) $  in this example in contradiction is  not trivial. \par

 \vspace{3mm}

Further, if an another distance $ r = r(x_1,x_2)  $ on the source set  $ X $ is such that

$$
\forall x_1 \in X \ \Rightarrow  \lim_{d(x,x_1) \to 0} \frac{d(x,x_1)}{r(x,x_1)} = 0, \eqno(3.4)
$$
then the space $ H^o(d) $ is continuously  embedded in the space $ H^o(r). $ \par
 We will write the equality (1.4) as follows: $  d << r. $\par

 For instance, the distance $  r(x_1,x_2) $ may has a form

$$
r(x_1, x_2) = d^{\beta}(x_1, x_2),  \ \beta = \const \in (0,1).
$$

 Of course, in the considered here problem $ X = [0, 2 \pi], \ d(\lambda, \mu) = |\lambda - \mu|^{\Delta}, \ \Delta = \const \in (0,1).  $
We introduce hence the following H\"older's spaces $ H^o_{\Delta}  $ over the circle $ [0, 2 \pi] $ consisting on all the periodical
(continuous) functions with finite norm

$$
||f||H^o_{\Delta} \stackrel{def}{=} = \sup_{\lambda \in [0, 2 \pi]} |f(\lambda)| +
\sup_{|\lambda -  \mu| > 0} \left\{ \frac{|f(\lambda) - f(\mu)|}{|\lambda - \mu|^{\Delta} }  \right\}, \eqno(3.5)
$$
and such that

$$
\lim_{\delta \to 0+} \left\{ \frac{\omega(f, h)}{h^{\Delta}} \right\} = 0.  \eqno(3.5a)
$$

 or equally
 $$
\forall \mu \in [0, 2 \pi] \ \Rightarrow
\lim_{\lambda \to \mu} \left\{ \frac{|f(\lambda) - f(\mu)|}{|\lambda - \mu|^{\Delta}} \right\} = 0.  \eqno(3.5b)
 $$

  The classical  CLT in H\"older's spaces,  i.e. CLT for the sums of independent random processes,
  with applications, is investigated in many works, see, e.g. \cite{Klicnarov’a1},
 \cite{Ostrovsky110}, \cite{Ratchkauskas1}, \cite{Ratchkauskas2}, \cite{Ratchkauskas5}. \par

 Our aim  in this section is investigation of the CLT for estimation of fractional derivative for spectral function. \par

\vspace{3mm}

{\bf Theorem 3.1.} {\it Let all the conditions of theorem 2.1 be satisfied.
Let also $ \Delta $ be arbitrary number such that $  0 \le \Delta < 1/2 - \alpha. $ The sequence of the distributions
generated in H\"older the space  $ H^o_{\Delta} $ by the r.p. $  \zeta_n $ converges weakly as $  n \to \infty $ to the
distribution in this space to at the same r.p. } $  \zeta_{\infty}. $ \par

\vspace{3mm}

{\bf Proof.}  To establish the weak compactness in these spaces, we return to  the inequalities (2.32)-(2.33):

$$
\omega(\zeta_n, h) \le C \ X_{n,p} \ \int_0^h \frac{t^{1/2 - 2/p} \ dt  }{t^{ 1 + \alpha  }} \le
$$

$$
C(\alpha,p) \ X_{n,p} \ h^{1/2 - \alpha - 2/p} \le C(\alpha,p) \ X_{n,p} \ h^{\Delta + \delta}, \  \exists
\delta = \const > 0, \eqno(3.6)
$$
if the value $ p  = \hat{p} =  \hat{p}(\alpha, \Delta, \delta)  $ is sufficiently great. As before, $ \sup_n {\bf E} X^p_{n,p} = 1.  $ \par

 We apply the Tchebychev's  inequality

$$
\sup_n {\bf P} \left( \frac{ \omega(\zeta_n, h) }{h^{\Delta + \delta}}  > u  \right) \le \frac{C(\alpha, \hat{p})}{ u^{\hat{p}} } \le \epsilon
\eqno(3.7)
$$
for sufficiently greatest values $  u. $  As long as the set of a (continuous) functions $  f: [0, 2 \pi] \to R $ such that

$$
\{f:  \omega(f,h) \le u \cdot h^{\Delta + \delta}  \}, \ u = \const < \infty
$$
is a shift - precompact set in the space $ H^o_{\Delta}, $  see \cite{Gusejnov1}, chapter 1, we conclude that the main
Prokhorov's condition  for weak  compactness of probability measures \cite{Prokhorov1} is satisfied. \par

 The rest: convergence of finite-dimensional distributions, belonging the limit process $ \zeta_{\infty}(\cdot) $ to the space
  $ H^o_{\Delta} $  is just.\par

\vspace{3mm}

{\bf Corollary 3.1.}  If we choose $ \alpha = 0, $ we get to the following extension of I.A.Ibragimov's \cite{Ibragimov1} result: for
arbitrary value  $  \Delta $ from the interval $ (0, 1/2) $ the sequence of distributions in the Banach space $ H^o_{\Delta} $
of the r.p. $  \tau_n(\cdot) $ converges weakly as $  n \to \infty $ to one for the r.p.  $ \tau(\cdot). $\par

\vspace{3mm}

{\bf Corollary 3.2.} If in addition to the conditions of theorem 3.1 the function $  F = F(\lambda) $ satisfies the following restriction

$$
\lim_{n \to \infty} \sqrt{n} \ ||  F^{(\alpha)} * \Phi_n - F^{(\alpha)}   ||H^o_{\Delta} = 0,  \eqno(3.8)
$$
then  the sequence of the distributions of H\"older continuous random processes

$$
  \theta_n(\lambda) \stackrel{def}{=} \sqrt{n} \ (F_{n,\alpha}(\lambda) - F^{(\alpha)}(\lambda))
$$
 converges weakly as $ n \to \infty $ in the space  $ H^o_{\Delta} $
 in the Prokhorov-Skorokhod sense to the distribution of at the same centered Gaussian continuous random process
$  \zeta_{\infty}(\cdot). $   \par

\vspace{3mm}

 The sufficient conditions  for the equality (3.8) may be found in the articles
\cite{Draganov1}, \cite{Lasuriya1}.  For instance, this equality is satisfied if

$$
D^{\alpha} F \in H^o_{\Delta},
$$
where as before $  \Delta < \alpha - 1/2,  $ see \cite{Draganov1}, page 8, corollary 3.2. \par

\vspace{4mm}

 \section{Non - asymptotical approach.}

\vspace{3mm}

 We do not suppose in this section that $  n \to \infty \ (n >> 1). $ More exactly, we intend to obtain here the upper and lower
 exponential estimate for the non-asymptotical probabilities  for the following normed uniform deviations

$$
W^o_{\alpha}(u) \stackrel{def}{=}  \sup_n {\bf P} ( \sqrt{n} \ \sup_{\lambda} | F_{n,\alpha}(\lambda)   - {\bf E} F_{n,\alpha}(\lambda) | > u ))
\eqno(4.1a)
$$
and  correspondingly
$$
W_{\alpha}(u) \stackrel{def}{=}   \sup_n {\bf P} ( \sqrt{n} \ \sup_{\lambda} | F_{n,\alpha}(\lambda)   - F^{(\alpha)}(\lambda) | > u )), \
u \ge 1. \eqno(4.1b)
$$

\vspace{3mm}

{\bf Theorem 4.1.} \ {\it Let all the conditions of theorem 2.1 be satisfied. Our statement: for some positive finite constants }
$  C_1 = C_1(\alpha), \ C_2 = C_2(\alpha), \ C_1 \le C_2 $
$$
\exp(- C_2(\alpha) u) \le  W^o_{\alpha}(u) \le \exp(- C_1(\alpha) u), \ u \ge 1. \eqno(4.2)
$$

\vspace{3mm}

{\bf Theorem 4.2.} \ {\it Let all the conditions of theorem 2.2 be satisfied. Our statement: for some positive finite constants }
$  C_3 = C_3(\alpha), \ C_4 = C_4(\alpha), \ C_3 \le C_4 $

$$
\exp(- C_4(\alpha) u) \le  W_{\alpha}(u) \le \exp(- C_3(\alpha) u), \ u \ge 1. \eqno(4.3)
$$

\vspace{3mm}

{\bf Proof.}  Let us consider the random processes $  \zeta_n = \zeta_n(\lambda). $ We employ the inequality (2.33):

$$
\omega(\zeta_n, h) \le C \ X_{n,p} \ \int_0^h \frac{t^{1/2 - 2/p} \ dt  }{t^{ 1 + \alpha  }},
$$
then

$$
\omega(\zeta_n, h) \le C_5 \ X_{n,p} \  h^{1/2 - 2/p}, \ p \ge 1/\alpha,
$$
where as before $ \sup_n \sup_p {\bf E} |X_{n,p}|^p = 1. $ On the other words,

$$
|\zeta_n(\lambda) - \zeta_n(\mu)| \le C_5 \ X_{n,p} \  |\lambda - \mu|^{1/2 - 2/p}. \eqno(4.4)
$$

 Analogously

$$
|\zeta_n(\lambda) | \le C_6 \ X_{n,p}. \eqno(4.4a)
$$

 As long as the random variables $ \{ \zeta_n(\lambda) \}, \ \lambda \in [0, 2 \pi]  $ are also the two-dimensional
stochastic integrals over Gaussian measure (Gaussian chaos),  on the other words, belongs to the described above Banach
space $  B(\phi). $  The so-called entropy condition \cite{Ostrovsky111} for the set $ [0, 2 \pi]  $ relative the distance
$  |\lambda - \mu|^{1/2 - 2/p} $ for each the values $  p, \ p > 1/\alpha  $ is satisfied, and we conclude using the main result
of an article \cite{Ostrovsky111} that

$$
\sup_n \max_{\lambda \in [0, 2 \pi] } |\zeta_n(\lambda)| \in G\psi, \eqno(4.5)
$$
where (recall) $ \psi(p) = p. $ \par
 The right-hand side of bilateral inequality (4.2) follows immediately from (4.5), see  \cite{Kozachenko1}, \cite{Ostrovsky1}.
The left-hand of (4.2) estimate  is very simple:

$$
 W^o_{\alpha}(u)  \ge {\bf P} ( |\zeta_1(\pi)| > u ) \ge \exp(- C_2(\alpha) u).
$$
 The second theorem 4.2 may be proved by means of theorem 2.2. \par

\vspace{4mm}

\section{Concluding remarks.}

\vspace{3mm}

{\bf A. Weight case. } \\

 Perhaps, it is interest to  investigate the error in the uniform norm of the approximation of a form

 $$
 V(x) \cdot ( W \cdot F)^{(\alpha)}(x) \approx V(x) \cdot ( W \cdot F)_{n,\alpha}(x),
 $$
or analogously

$$
 V(x) \cdot ( W * F)^{(\alpha)}(x) \approx V(x) \cdot ( W * F)_{n,\alpha}(x),
$$
or analogously

$$
 V(x) * ( W * F)^{(\alpha)}(x) \approx V(x) *( W * F)_{n,\alpha}(x),
$$
where $ V(x), W(x)  $ are two weight functions, for instance, $ V(x) = |x|^{\gamma}, \ W(x) = |x|^{\Delta},  \hspace{6mm}
\gamma, \Delta = \const.  $ \par

\vspace{3mm}

{\bf B. Applications (possible) in  statistics.}\\

 The asymptotical  tail behavior of the  statistic $ \sup_{\lambda} |\theta_{n}(\lambda)| $ may be used perhaps  in turn in statistics,
for instance, for the verification of semi-parametrical hypotheses and detection of distortion times of signals etc. \par

\vspace{3mm}

{\bf C.  Non - centered sample.} \\

 If the source stationary Gaussian random sequence $  \{  \eta_k  \}, \ k = 1,2,\ldots $ is non - centered:

$$
{\bf E} \eta_k = a \ne 0, \ k = 1,2, \ldots,n,
$$
then we can replace as ordinary

$$
\eta_k := \eta_k^o \stackrel{def}{=} \eta_k - n^{-1}\sum_{j=1}^n \eta_j = \eta_k - a_n,
$$
where $ a_n = n^{-1}\sum_{j=1}^n \eta_j  $ is consistent estimation for the value $  a. $ Both the theorems 2.1 and 2.2 remains true
under at the same conditions.\par

\vspace{3mm}

{\bf D.} Perhaps, obtained above results may be extended on the multivariate time series by using of the results of the book
\cite{Hannan1}, chapter 5, as well as on the non-Gaussian processes through the cumulant function
 and on the case of the "continuous time". \par

\end{document}